\documentclass[12pt,leqno]{article}
\usepackage{amsfonts}
\usepackage{amsmath, amsthm, amsfonts, amssymb, color}
\usepackage{mathrsfs}
\usepackage{color}
\usepackage[notcite,notref]{showkeys}

\theoremstyle{definition}

\newcommand{\scr}[1]{\mathscr #1}
\definecolor{wco}{rgb}{0.5,0.2,0.3}

\numberwithin{equation}{section} \theoremstyle{remark}

\newcommand{\ua}{\uparrow}

\title{{\bf  Singular McKean-Vlasov  SDEs:  Well-Posedness,  Regularities and Wang's Harnack Inequality}\footnote{Supported in
 part by  ERC-AdG 694405 "RicciBounds".} }
\author{
{\bf Panpan Ren
 }\\
\footnotesize{ Mathematics department, Bonn University,  Germany}\\
\footnotesize{ rppzoe@gmail.com, 
}}
\begin{document}
\allowdisplaybreaks
\def\R{\mathbb R}  \def\ff{\frac} \def\ss{\sqrt} \def\B{\mathbf
B}
\def\N{\mathbb N} \def\kk{\kappa} \def\m{{\bf m}}
\def\ee{\varepsilon}\def\ddd{D^*}
\def\dd{\delta} \def\DD{\Delta} \def\vv{\varepsilon} \def\rr{\rho}
\def\<{\langle} \def\>{\rangle}
  \def\nn{\nabla} \def\pp{\partial} \def\E{\mathbb E}
\def\d{\text{\rm{d}}} \def\bb{\beta} \def\aa{\alpha} \def\D{\scr D}
  \def\si{\sigma} \def\ess{\text{\rm{ess}}}\def\s{{\bf s}}
\def\beg{\begin} \def\beq{\begin{equation}}  \def\F{\scr F}
\def\Ric{\mathcal Ric} \def\Hess{\text{\rm{Hess}}}
\def\e{\text{\rm{e}}} \def\ua{\underline a} \def\OO{\Omega}  \def\oo{\omega}
 \def\tt{\tilde}\def\[{\lfloor} \def\]{\rfloor}
\def\cut{\text{\rm{cut}}} \def\P{\mathbb P} \def\ifn{I_n(f^{\bigotimes n})}
\def\C{\scr C}      \def\aaa{\mathbf{r}}     \def\r{r}
\def\gap{\text{\rm{gap}}} \def\prr{\pi_{{\bf m},\varrho}}  \def\r{\mathbf r}
\def\Z{\mathbb Z} \def\vrr{\varrho} \def\ll{\lambda}
\def\L{\scr L}\def\Tt{\tt} \def\TT{\tt}\def\II{\mathbb I}
\def\i{{\rm in}}\def\Sect{{\rm Sect}}  \def\H{\mathbb H}
\def\M{\mathbb M}\def\Q{\mathbb Q} \def\texto{\text{o}} \def\LL{\Lambda}
\def\Rank{{\rm Rank}} \def\B{\scr B} \def\i{{\rm i}} \def\HR{\hat{\R}^d}
\def\to{\rightarrow}\def\l{\ell}\def\iint{\int}\def\gg{\gamma}
\def\EE{\scr E} \def\W{\mathbb W}
\def\A{\scr A} \def\Lip{{\rm Lip}}\def\S{\mathbb S}
\def\BB{\scr B}\def\Ent{{\rm Ent}} \def\i{{\rm i}}\def\itparallel{{\it\parallel}}
\def\g{{\mathbf g}}\def\Sect{{\mathcal Sec}}\def\T{\mathcal T}\def\BB{{\bf B}}
\def\f{\mathbf f} \def\g{\mathbf g}\def\BL{{\bf L}}  \def\BG{{\mathbb G}}
\def\Bd{{D^E}} \def\BdP{D^E_\phi} \def\Bdd{{\bf \dd}} \def\Bs{{\bf s}} \def\GA{\scr A}
\def\Bg{{\bf g}}  \def\Bdd{\psi_B} \def\supp{{\rm supp}}\def\div{{\rm div}}
\def\ddiv{{\rm div}}\def\osc{{\bf osc}}\def\1{{\bf 1}}\def\BD{\mathbb D}\def\GG{\Gamma}
\def\H{{\bf H}}
\maketitle

\begin{abstract}  The well-posedness and regularity  estimates  in initial distributions are derived   for   singular McKean-Vlasov SDEs, where the drift contains a locally standard integrable term and a superlinear term in the spatial variable, and is Lipchitz continuous in the distribution variable with respect to a weighted variation distance. When the superlinear term is strengthened to be Lipschitz continuous, Wang's Harnack inequality is established.
These results   are new also for the classical It\^o SDEs where the coefficients are distribution independent.     \end{abstract} \noindent
 AMS subject Classification:\  60B05, 60B10.   \\
\noindent
 Keywords:  Singular,  McKean-Vlasov SDE,  Regular,  Wang's Harnack inequality.

 \vskip 2cm

 \section{Introduction and main results}
 
  In recent years, singular SDEs have been intensively investigated by using   Zvonkin's transform \cite{ZV} and  Krykov's estimate \cite{Krylov}   developed from \cite{Ver} for bounded drift, \cite{KR,Z2} for integrable drift,
 and \cite{XXZZ, YZ0} for locally integrable drift.   
 
 In this paper, we aim to improve existing results on the well-posedness and regularity estimates derived for singular SDEs,  and make extensions to McKean-Vlasov SDEs
 (also called mean field SDEs or distribution dependent SDEs),
a hot research object due to its essential links to nonlinear Fokker-Planck equations and mean field particle systems, see for instance the monographs \cite{SN}, \cite{CD} and  the survey \cite{20HRW}.

  Let  $\scr  P$ be the space of probability measures on $\R^d$ equipped with the weak topology. 
For $T\in (0,\infty)$, consider 
 \beq\label{E10} \d X_t= b_t(X_t,\L_{X_t})\d t+\si_t(X_t)\d W_t,\ \ t\in [0,T],\end{equation}
where $W_t$ is an $m$-dimensional Brownian motion on a complete filtration probability space $(\OO, \{\F_t\}_{t\ge 0},\P)$,   $\L_\xi$ is the distribution (i.e. the law) of a random variable $\xi$, and 
\beg{align*}b: [0,\infty)\times\R^d\times \scr  P\to \R^d, 
 \ \ \si: [0,\infty)\times\R^d \to \R^d\otimes\R^m\end{align*}  are measurable. 

To characterize the dependence on the distribution variable, we introduce some probability distances. 
For a measurable function $V\ge 1$, 
let $$\scr P_V:=\bigg\{\mu\in \scr P: \mu(V):=\int_{\R^d} V\d\mu<\infty\bigg\}$$ be equipped with the $V$-weighted variation metric
$$\|\mu-\nu\|_V:= \sup_{|f|\le V} |\mu(f)-\nu(f)|,\ \ \mu,\nu\in \scr P_V.$$
When $V=1$ it reduces to the total variation norm, i.e. $\|\cdot\|_{1}= \|\cdot\|_{var}$.   We often take $V$ to be a compact function, i.e. its level sets $\{V\le r\}$ for $r>0$ are compact. 

Next, for any $k\in [1,\infty)$, let 
$\scr P_k=\scr P_V$ for $V:= 1+|\cdot|^k$, i.e.
$$\scr P_k:=\{\mu\in \scr P: \|\mu\|_k:= \mu(|\cdot|^k)^{\ff 1 k} <\infty\}.$$  In this case, we denote $\|\cdot\|_V=\|\cdot\|_{k,var}$. 
Consider the  $L^k$-Wasserstein distance 
$$\W_k(\mu,\nu):=   \inf_{\pi\in \C(\mu,\nu)} \bigg(\int_{\R^d\times\R^d} |x-y|^k \pi(\d x,\d y) \bigg)^{\ff 1{k}},\ \ \mu,\nu\in \scr P_ k,$$ 
where $\C(\mu,\nu)$ is the set of all couplings for $\mu$ and $\nu$.   According to \cite[Theorem 6.15]{VV},   there exists a constant $c>0$ such that
$$ \|\mu-\nu\|_{var}+\W_k(\mu,\nu)^{k} \le c \|\mu-\nu\|_{k,var},\ \ \mu,\nu\in \scr P_k.$$ 

To measure the singularity of the drift, we recall some functional spaces introduced in \cite{XXZZ}.  For any $p\ge 1$, $L^p(\R^d)$ is the class of   measurable  functions $f$ on $\R^d$ such that \index{$L^p(\R^d)$: $L^p$-space for Lebesgue measure on $\R^d$}
 $$\|f\|_{L^p(\R^d)}:=\bigg(\int_{\R^d}|f(x)|^p\d x\bigg)^{\ff 1 p}<\infty.$$
For any $\epsilon >0$ and $p\ge 1$, let  $H^{\epsilon,p}(\R^d):=(1-\DD)^{-\ff\epsilon 2} L^p(\R^d)$ with
 $$\|f\|_{H^{\epsilon,p}(\R^d)}:= \|(1-\DD)^{\ff\epsilon 2} f\|_{L^p(\R^d)}<\infty,\ \ f\in  H^{\epsilon,p}(\R^d).$$ \index{$H^{\epsilon,p}(\R^d)$: domain of $(1-\DD)^{\ff\epsilon 2} $ in $L^p(\R^d)$}
  
  For any $z\in\R^d$ and $r>0$,   let $B(z,r):=\{x\in\R^d: |x-z|< r\}$ be the open ball centered at $z$ with radius $r$.
For any $p,q\geq1$,  let $\tt L_q^p$ denote the class of measurable functions $f$ on $[0, T]\times\R^d$ such that
$$\|f\|_{\tt L_q^p}:= \sup_{z\in \R^d}\bigg( \int_{0}^{T} \|1_{B(z,1)}f_t\|_{L^p(\R^d)}^q\d t\bigg)^{\ff 1 q}<\infty.$$
For any $\epsilon>0$, let $\tt H_{q}^{\epsilon,p}$ be the space of $f\in \tt L_q^p$ with
$$\|f\|_{\tt H_q^{\epsilon,p}}:=  \sup_{z\in \R^d}\bigg( \int_{0}^{T}\|g(z+\cdot)f_t\|_{H^{\epsilon,p} (\R^d)}^q\d t\bigg)^{\ff 1 q}<\infty$$ for some $g\in C_0^\infty(\R^d)$
satisfying $g|_{B(0,1)}=1,$ where $C_0^\infty(\R^d)$ is the class of $C^\infty$ functions on $\R^d$ with compact support.
We remark that the space $\tt H_{q}^{\epsilon,p}$ does not depend on the choice of $g$.

We will  take $(p,q)$ from the class
$$\scr K:=\Big\{(p,q): p,q\in  (2,\infty),\  \ff d p+\ff 2 q<1\Big\}.$$

\subsection{Well-posedness} 

Let us first recall the definition of well-posedness.

\beg{defn} A continuous adapted process  $(X_{s,t})_{T\ge t\ge s}$   is called a  solution of \eqref{E10} from time   $s$, if
$$\int_s^t\E\big[|b_r(X_{s,r}, \L_{X_{s,r}})|+  \| \si_r(X_{s,r})\|^2\big]\d r<\infty,\ \ T\ge t\ge s,$$ and   $\P$-a.s.
  $$X_{s,t} = X_{s,s} +\int_s^t b_r(X_{s,r}, \L_{X_{s,r}})\d r + \int_s^t \si_r(X_{s,r})\d W_r,\ \ T\ge t\ge s.$$
When $s=0$ we simply denote   $X_t=X_{0,t}$.

 A couple $(\tt X_{s,t},\tt W_t)_{T\ge t\ge s}$   is called a  weak solution of \eqref{E10} from time $s$, if   $\tt W_t$  is the   $m$-dimensional Brownian motion on a complete filtration probability space $(\tt\OO,\{\tt\F_t\}_{t\in [0,T]}, \tt\P)$  such that   $(\tt X_{s,t})_{T\ge t\ge s}$  is a solution of \eqref{E10}  from time   $s$
for   $(\tt W_t, \tt\P)$ replacing   $(W_t,\P)$.    \eqref{E10} is called weakly unique for an initial distribution $\nu\in\scr  P$, if all weak solutions with distribution $\nu$ at time $s$ are equal in law.

Let $\hat {\scr  P}$ be a subspace of $\scr  P$. \eqref{E10} is called  strongly (respectively,  weakly) well-posed for   distributions in $\hat {\scr  P}$, if for any  $s\in [0,T)$ and $\F_s$-measurable $X_{s,s}$  with $\L_{X_{s,s}}\in \hat {\scr  P}$ (respectively, any initial distribution $\nu\in \hat {\scr  P}$ at time $s$), it has a unique strong (respectively, weak) solution.
  
    \eqref{E10}  is called well-posed  for distributions in $\hat {\scr P}$,  if it is both  strongly and weakly well-posed for distributions in $\hat {\scr P}$.
\end{defn}

To prove the well-posedness, we make the following assumptions.

 \beg{enumerate}\item[$(H_1)$]  $\si_t(x)$ and $b_t(x,\mu)=b^{(0)}_t(x)+ b_t^{(1)}(x,\mu)$ satisfy the following conditions for  a compact function $1\le V\in C^2(\R^d; [1,\infty))$. 
 \item[$(1)$]  $a:= \si\si^*$ is invertible with $\|a\|_\infty+\|a^{-1}\|_\infty<\infty$, where $\si^*$ is the transposition of $\si$,   and
$$ \lim_{\vv\to 0} \sup_{|x-y|\le \vv, t\in [0,T]} \|a_t(x)-a_t(y)\|=0.$$

 \item[$(2)$]   $|b^{(0)}|\in \tt L_{q_0}^{p_0}$ for some  $(p_0,q_0)\in \scr K$,  Moreover, $\si_t$ is weakly differentiable such that 
\beq\label{ASI} \|\nn \si_t\|\le \sum_{i=1}^l f_i \end{equation}  holds for some $l\in \mathbb N$ and $0\le f_i \in \tt L_{q_i}^{p_i}$ with $  (p_i,q_i) \in  \scr K, 1\le i\le l.$  
  \item[$(3)$]   for any $\mu\in C([0,T];\scr P_V)$,  $b^{(1)}_t(x,\mu_t)$   is  locally bounded in $(t,x)\in [0,T]\times\R^d$. Moreover,  there exist constants $K,\vv>0$ a compact function $V\in C^2(\R^d; [1,\infty))$     such that
\beg{align*} &\sup_{|y-x|\le \vv}\big\{|\nn V(y)|+   \|\nn^2 V(y)\| \big\}\le K V(x),\\
&   \<b_t^{(1)}(x,\mu),\nn V(x) \>+\vv |b_t^{(1)}(x,\mu)| \sup_{B(x,\vv)} \big\{|\nn V|+ |\nn^2 V \| \big\} \\
& \le K \big\{V(x)+ \mu(V)\big\},\ \ \ 
   x\in \R^d,\mu\in \scr P_V. \end{align*}  
\item[$(4)$]   there exists a constant $\kk>0$ such that 
 \beq\label{KK1}  |b_t(x,\mu)- b_t(x,\nu)|\le \kk \|\mu-\nu\|_V,\ \ \mu,\nu\in \scr P_V, x\in\R^d.\end{equation} 
\end{enumerate}

\beg{thm}\label{TA1}  Assume $(H_1)$.   Then $\eqref{E10}$  is well-posed for distributions in $\scr P_V$. Moreover: \beg{enumerate} 
\item [$(1)$]   for any $n\ge 1$  there exists a constant $c(n)>0$ such that     
 \beq\label{KK2} \E\Big[\sup_{t\in [0,T]} V(X_t)^n\Big|X_0\Big]\le c(n) \big\{(\E[ V(X_0)] )^n + V(X_0)^n\big\}\end{equation}  holds  for any solution $X_t$ of $\eqref{E10}$ with $\L_{X_0}\in \scr P_V$.
\item [$(2)$]    for any sequence $\{\mu_n\}_{n\ge 1}\subset \scr P_V$ with bounded $\mu_n(V^p)$ for some $p>1$ such that $\mu_n\to\mu$ weakly, 
 \beq\label{KK3} \lim_{n\to\infty} \|P_t^*\mu_n- P_t^*\mu\|_V=0.\end{equation}  
 \item[$(3)$] if there exists a constant $K>0$ such that 
\beq\label{KK1'} |b_t(x,\mu)- b_t(x,\nu)|\le K\|\mu-\nu\|_{var}\ \ \mu,\nu\in \scr P_V,\end{equation}  
then    $\eqref{KK2}$ holds for some constant $c>0$ independent of $\mu$, and 
 \beq\label{KK3'} \lim_{\nu\to\mu \ \text{weakly}} \|P_t^*\mu- P_t^*\nu\|_{var}=0.\end{equation} \end{enumerate} 
\end{thm}

\paragraph{Remark 1.1.}    (1) Theorem \ref{TA1} extends existing  well-posedness results derived for singular McKean-Vlasov SDEs, for instance:
   \beg{enumerate}\item[$(a)$]  \cite{HW20} and \cite{HW20c}  with $l=1$  in $(H_1)$(2), and 
    $(H_1)$(4)  with $V:= (1+|\cdot|^2)^{\ff k 2}$,  and the following stronger condition  stronger than $(H_1)(3)$: 
 \beq\label{KK'}  \sup_{t\in [0,T], x\ne y  } \bigg\{|b_t^{(1)}(0,\dd_0)|+ \ff{|b_t^{(1)}(x,\dd_0)-b_t^{(1)}(y,\dd_0)|}{|x-y|}\bigg\}<\infty,\end{equation}
 where $\dd_0$ is the Dirac measure at $0$. 
\item[$(b)$]    \cite{RZ} as well as \cite{Zhao} for  H\"older continuous $\si_t$ and $\sup_\mu \|b(\cdot,\mu)\|_{\tt L_{q_0}^{p_0}}<\infty$. 
\end{enumerate} 
 
 (2) The regularity property  included in \eqref{KK3} and \eqref{KK3'}  is  new in this general situation. Under \eqref{KK'} replacing $(H_1)$(4), 
  the log-Harnack inequality was established in \cite{W21b} so that 
 $$ \|P_t^*\mu- P_t^*\nu\|_{var}\le \ff c {\ss t} \W_2(\mu,\nu)$$ holds for some constant $c>0$, which is incomparable with \eqref{KK3} since $\|\cdot\|_{var}$ is essentially smaller than $\|\cdot\|_{V}$. 
  
  (3)  Theorem \ref{TA1} is new even in the setting of singular SDEs, see comments before Theorem \ref{T01}.



\subsection{Wang's  Harnack inequality}

Since 1997 when the dimension-free Harnack inequality  of type 
$$|P_t f(x)|^p\le (P_t |f|^p(y)) \e^{\ff c t \rr(x,y)^2}$$ was  found in   Wang \cite{W97} for diffusion semigroups $P_t$ on Riemannian manifolds, 
this type inequality  has been intensively developed and applied to many different models, see \cite{W13} 
for a general theory on the study.  In recent years, Wang's   inequality has been established for McKean-Vlasov SDEs in \cite{W18} under monotone conditions  as well as  in \cite{HW19} for bounded $b$
which is Dini continuous in the space variable and $\W_2$-Lispchitz continuous in the distribution variable. 

In this paper,  we establish dimension-free  Harnack inequality in a more general situation, which is new even 
for classical SDEs, see comments before Theorem \ref{T02} in the next section. 

 \beg{enumerate}\item[$(H_2)$]   $(H_1)$(1)-(2),   $\eqref{KK'}$ and the following conditions hold. \item[$(1)$] 
  There exists increasing $\Phi\in C^2([0,\infty);  [1,\infty))$ with 
\beq\label{PHII} \limsup_{r\to\infty} \ff{\Phi'(r)+|\Phi''(r)|}{\Phi(r)}<\infty,\end{equation}
such that for some constant $\kk>0$ and  $V:= \Phi(|\cdot|^2), $  
\beq\label{MOO} |b_t(x,\mu)-b_t(x,\nu)| \le \kk \|\mu-\nu\|_V,\ \ (t,x,\mu)\in [0,T]\times \R^d\times \scr P_V.\end{equation} 
\item[$(2)$] There exists increasing $\varphi\in C([0,\infty); [0,\infty))$ satisfying  $\varphi(0)=0, \varphi(r)>0$ for $r>0$,   $\psi(r):= \ff {r^2} {\varphi(r)^2}$ is increasing in $r>0$ and 
$\int_0^1 \ff{(\varphi\circ\psi^{-1})^2(s)} s\d s<\infty, $ such that 
$$\|\si_t(x)-\si_t(y)\|\le \varphi(|x-y|),\ \ x,y\in \R^d, t\in [0,T].$$
  \end{enumerate} 
  
Typical examples of $\varphi$ in $(H_2)$(2) include $\varphi(r)=r^\aa$ for $\aa \in (0,1)$  and 
 $\varphi(r)=\log^{-\theta}(\e +r^{-1})$ for $\theta> 1$, where in the first case $\si_t$ is H\"older continuous and in the second case it is only Dini continuous.

 By Theorem \ref{TA1}, $(H_2)$ implies the well-posedness of $\eqref{E10}$ for distributions in $\scr P_V$. Consider 
 $$P_t f(\mu):= \int_{\R^d} f\d(P_t^*\mu),\ \ t\in [0,T], f\in \B_b(\R^d), \mu\in \scr P_V.$$
 
\beg{thm}\label{TA2} Assume $(H_2)$.  Then the following assertions hold. 
\beg{enumerate} \item[$(1)$] There exist constants $c,p>1$ such that for any $ \ t\in (0,T], f\in \B_b(\R^d),$  
\beq\label{HHH} |P_tf|^p(\mu)\le\{ P_t |f|^p(\nu)\} \inf_{\pi\in \C(\mu,\nu)} \int_{\R^d\times\R^d}\e^{c+\ff c t |x-y|^2}\pi(\d x,\d y),\ \ \mu,\nu\in \scr P_V.\end{equation} 
\item[$(2)$] If $\Phi$ is bounded then there exists a constant $c>0$ such that for any $t\in (0,T],$ 
\beq\label{GRR} \|P_t^*\mu-P_t^*\nu\|_{var}^2 \le c\big(t^{-1}-\log [1\land \W_2(\mu,\nu)]\big) \W_2(\mu,\nu)^2,\ \  \mu,\nu\in \scr P_V.\end{equation} \end{enumerate} 
\end{thm} 

\paragraph{Remark 1.2.} (1) By the proof of  \cite[Theorem 1.4.2]{W13}, if the right hand side in \eqref{HHH} is finite, then  $\rr_t^{\mu,\nu}:=\ff{\d P_t^*\mu}{\d P_t^*\nu}$ exists and satisfies 
$$\big\{P_t (\rr_t^{\mu,\nu})^{\ff 1 {p-1}}(\mu)\big\}^{p-1} \le  \inf_{\pi\in \C(\mu,\nu)}\int_{\R^d\times\R^d} \e^{c+\ff c t |x-y|^2}\pi(\d x,\d y).$$
The Harnack inequality \eqref{HHH} is new even for the classical distribution dependent SDEs, see comments before Theorem \ref{T02}. 

(2) By \eqref{GRR}, for any $f\in\B_b(\R^d)$ and $t\in (0,T]$, $P_tf$ is nearly Lipschitz continuous in $\W_2$ in the sense that
$$|P_tf(\mu)-P_tf(\nu)|\le \|f\|_\infty \ss{c\big(t^{-1}-\log [1\land \W_2(\mu,\nu)]\big)}\, \W_2(\mu,\nu).$$
When \eqref{MOO} holds for $\W_2(\mu,\nu)$ replacing $\|\mu-\nu\|_V$,   Theorem 4.1 in \cite{W21b} implies the exact Lipschitz continuity of $P_tf$ in $\W_2$. 
See also \cite{HW19} and   \cite{HW21}  for the $\W_2$-Lipschitz continuity of $P_tf$  under stronger conditions on $b$, where \cite{HW21} allows $\si$ to be distribution dependent. 

\  

In the following three  sections,  we first prove the above results for singular SDEs where $b_t(x,\mu)=b_t(x)$ does not depend on $\mu$, then extend to the distribution dependent setting to prove the above two theorems.  

 \section{Singular SDEs}

  Consider the following  SDE on $\R^d$:
 \beq\label{E1} \d X_t= b_t(X_t)\d t +\si_t(X_t)\d W_t,\ \ t\in [0,T].\end{equation}
  There are a plenty of papers studying the well-posedness of this SDE. In the following  we mention two typical results  under   weak  monotone  condition and locally integrable condition respectively. 

 According to  \cite{FZ},  when $b$ and $\si$ are  continuous  satisfying the following  weak semi-Lipschitz continuous condition:
\beq\label{SMM} \beg{split} &2\<x-y,  b_t(x)-b_t(y)\>+\|\si_t(x)-\si_t(y)\|_{HS} ^2\\
&\le K |x-y|^2\log (2+|x-y|^{-1}),\\
&2\<b_t(x), x\>+\|\si_t(x)\|_{HS}^2 \le K(1+|x|^2)\log (2+|x|^2),\ \   \ t\in [0,T], x,y\in \R^d,\end{split}\end{equation}  then \eqref{E1}  is well-posed.

On the other hand, in recent years \eqref{E1} has been intensively studied under locally integrable conditions.  According to \cite{YZ0}, see \cite{XXZZ,Z2} 
and references within for earlier results, the well-posedness of \eqref{E1}  holds under the following assumption. We remark that these papers (also related existing    references) only consider 
the case $l=1$ in condition (2) below, but the proof applies to   $l\ge 2$ by replacing  $\|\nn \si\|$ with $\sum_{i=1}^l f_i$ and applying Khasminskii's estimate to each   $f_i$ respectively.

\beg{enumerate}
\item[$(A_1)$]           Let    $a_t(x):= (\si_t\si^*_t)(x)$ and $b_t(x)= b_t^{(0)}(x)+ b_t^{(1)}(x)$.
\item[$(1)$] $a$  is invertible with $\|a\|_\infty+\|a^{-1}\|_\infty<\infty$    and uniformly continuous in $x$: 
$$ \lim_{\vv\to 0} \sup_{|x-y|\le \vv, t\in [0,T]} \|a_t(x)-a_t(y)\|=0.$$
\item[$(2)$] There exist $l\in \mathbb N$,  $\{(p_i,q_i)\}_{0\le i\le l} \subset \scr K$  and $0\le  f_i \in \tt L_{q_i}^{p_i}, 1\le i\le l$    such that
$$  |b^{(0)}|\in \tt L_{q_0}^{p_0},\ \ \|\nn\si\| \le \sum_{i=1}^l f_i .$$
 \item[$(3)$]   $b_t^{(1)}$ is Lipschitz continuous with  
 $$ \sup_{t\in [0,T]}\big\{|b^{(1)}_t(0)|+\|\nn b_t^{(1)}\|_\infty\big\}<\infty,$$   where   $\|\nn b_t^{(1)}\|_\infty$ is  the Lipschitz constant of   $b_t^{(1)}$.\end{enumerate}

\paragraph{Remark 2.1.} $(A_1)$ does not include \eqref{SMM}. Our first result ensures  the well-posedness  under the following condition,
which extends both  $(A_1)$  and \eqref{SMM}. Indeed, when  $\si$ is bounded, $(3')$ holds for $V(x)=1+|x|^2$ if there exist  constants $\vv,C>0$ such that  
 $$\<b_t^{(1)}(x),x\>+ \vv |b_t^{(1)}(x)| \le C \phi(1+|x|^2),\ \ x\in\R^d, t\in [0,T],$$
 which extends  \eqref{SMM} by allowing a singular term for $\phi(s):= \log (\e +s)$ and covers  $(A_1)$  for $\phi(s):=1+s$.

\beg{enumerate}
\item[$(A_2)$]        Assume $(A_1)(1)$-$(2)$  and $b^{(1)}$ is locally bounded such that the following condition holds. 
\item[$(3')$] There exist constants $K,\vv>0$, increasing $\phi\in C^1([0,\infty); [1,\infty))$ with $\int_0^\infty\ff{\d s}{s+\phi(s)}=\infty$, 
and   $V\in C^2(\R^d; [1,\infty))$  having compact level sets,  such that 
\beg{align*} &\sup_{|y-x|\le \vv}\big\{|\nn V(y)|+   \|\nn^2 V(y)\| \big\}\le K V(x),\\
&   \<b^{(1)}(x),\nn V(x) \>+\vv |b^{(1)}(x)| \sup_{B(x,\vv)} \Big\{ |\nn V|+ |\nn^2 V|  \Big\}\le K\phi(V(x)),  \ \ x\in \R^d.\end{align*} \end{enumerate}

 In the following, we prove the well-posedness with strong Feller property and moment estimates  under assumption $(A_2)$ and establish Wang's Harnack inequality under $(A_1)$. 
 
 \subsection{Well-posedness} 
  
 \beg{thm} \label{T01} Assume $(A_2)$.  
 Then $\eqref{E1} $ is well-posed.  Moreover:
 \beg{enumerate}
 \item[$(1)$]  For any $n\ge 1$ and $B_n:=\{x\in \R^d: |x|\le n\}$, 
$$ \lim_{\vv\downarrow 0} \sup_{x,y\in B_n, |x-y|\le\vv } \E\Big[\sup_{t\in [0,T]} |X_t^x-X_t^y|\land 1\Big]=0,$$
 where $X_t^x$ is the solution starting at $x$. 
 \item[$(2)$]   Let $P_t^*\dd_x=\L_{X_t^x}$ be the distribution of $X_t^x$.  Then
 $$\lim_{y\to x} \|P_t^*\dd_x-P_t^*\dd_y\|_{var}=0,\ \ t\in (0,T], x\in \R^d.$$ Equivalently, the associated semigroup $(P_t)_{t\in (0,T]}$ is strong Feller, i.e. $P_t \B_b(\R^d)\subset C_b(\R^d)$. 
  \item[$(3)$] If $\phi(r)=r$, then for any $k\ge 1$ there exists a constant $c(k)>0$ such that 
 $$\E\Big[\sup_{t\in [0,T]} V(X_t^x)^k \Big]\le c(k) V(x)^k,\ \ x\in \R^d.$$ 
\end{enumerate}  \end{thm} 
 
\beg{proof} 
     (a) For any $n\ge 1$, let  $B_n:=\{x: |x|\le n\}$ and 
$$b^n_t:= 1_{B_n}b_t^{(1)}+b_t^{(0)},\ \ t\in [0,T].$$
By Theorem 1.1 in \cite{XXZZ}, for any $x\in \R^d$, the following SDE is well-posed:
\beq\label{E1'} \d X_t^n= b^n(X_t^n)\d t +\si(X_t^n)\d W_t,\end{equation}
and for $X_t^{x,n}$ being the solution starting at $x$, 
\beq\label{C1} \sup_{x\ne y} \E\Big[\sup_{t\in [0,T]} \ff{|X_t^{x,n}-X_t^{y,n}|^k}{|x-y|^k} \Big]<\infty,\ \ k\ge 1.\end{equation} 
 As we already mentioned before that \cite{XXZZ} only considers $l=1$ in condition $(A_1)$(2),  but the proof  works also for $l\ge 2$ by applying Khasminskii's estimate to
 each $f_i$  replacing $\|\nn\si\|$. 

Let $\tau_n^x:= \inf\{t\ge 0: T\land |X_t^{x,n}|\ge n\}.$ Then $X_t^{x,n}$ solves \eqref{E1} up to time $\tau_n^x$, and by the uniqueness we have 
$$X_t^{x,n}=X_t^{x,m},\ \ t\le \tau_n^x\land\tau_m^x, n,m\ge 1.$$
So, it suffices to prove that $\tau_n^x\to T$ as $n\to\infty$. 

Let $$L_t := \ff 1 2{\rm tr}\big\{\si_t\si_t^* \nn^2\big\}+\nn_{b_t^{(0)}}.$$ By \cite[Theorem 3.1]{XXZZ} and $(A_1)(1)$-$(2)$, for any $\ll\ge 0$ the PDE
\beq\label{PDE0} (\pp_t+L_t )u_t=\ll u_t-b_t^{(0)},\ \ t\in [0,T], u_T=0\end{equation}
has a unique solution $u\in \tt H_{q_0}^{p_0} (T)$, and there exist constants $\ll_0, c,\theta>0$ such that 
\beq\label{ESS} \ll^{\theta} (\|u\|_\infty+\|\nn u\|_\infty) +\|\pp_t u\|_{\tt L_{q_0}^{p_0}(T)}+\|\nn^2 u\|_{\tt L_{q_0}^{p_0}(T)}\le c,\ \ \ll\ge \ll_0.\end{equation}
So, we may take $\ll\ge \ll_0$ such that 
\beq\label{*4}\|u\|_\infty+\|\nn u\|_\infty\le \vv.\end{equation} 
Let $\Theta_t(z)= z+u_t(z)$ for $(t,z)\in [0,T]\times\R^d$.    By It\^o's formula in \cite[Theorem 4.1(ii)]{XXZZ},  $Y_t^n:= \Theta(X_t^n)$ satisfies
\beq\label{YNN0}  \d Y_t^n = \big\{1_{B_n}b^{(1)} +\ll u_t+1_{B_n}\nn_{b^{(1)}} u_t  \big\}(X_t^n) \d t + \{(\nn \Theta_t)\si\}(X_t^n) \d W_t.\end{equation}
By \eqref{*4} and $(A_2)(3')$, there exist  $c_0,c_1,c_1>0$ such that for some martingale $M_t$, 
\beg{align*} &   \d  \{V(Y_t^n) +M_t\} \le \Big[ \big\<\{b^{(1)}+\nn_{b^{(1)}} u_t\}(X_t^n), \nn  V(Y_t^n)\big\>1_{B_n}+  c_0(|\nn V(Y_t^n)|+ \|\nn^2 V(Y_t^n)\|)\Big]\d t \\
&\le \Big\{\<b^{(1)}(X_t^n), \nn V(X_t^n)\> + \vv|b^{(1)}(X_t^n)| \sup_{B(X_t^n,\vv)} \Big(|\nn V|+   \|\nn^2 V)\| \Big)+c_0K V(Y_t^n)\Big)\Big\}\d t\\
&\le \big\{K \phi(V(X_t^n)) + c_0K V(Y_t^n)\big\}\d t \le K\big\{\phi((1+\vv K) V(Y_t^n)) + c_0 V(Y_t^n)\big\}\d t,\ \ t\le\tau_n^x.
  \end{align*}

 Let $H(r):=\int_0^r\ff{\d s}{s+ \phi((1+\vv K)s)}.$  Then $\int_0^\infty\ff{\d s}{s+\phi(s)}=\infty$ implies 
 \beq\label{HIN} H(\infty):=\lim_{r\to\infty} H(r)=\infty.\end{equation}
 Since $\phi\in C^1([0,\infty);[1,\infty))$ is increasing, by It\^o's formula we obtain
 $$\d H(V(Y_t^n))\le c_3\d t +\d \tt M_t,\ \ t\in [0,\tau_n^x]$$ for some constant $c_3>0$ and some martingale $\tt M_t$. 
 Then $$\E [(H\circ V)(Y_{t\land\tau_n^x}^n)]\le V(x+u_0(x)) +c_3 t,\ \ n\ge 1, t\in [0,T].$$
Since  \eqref{*4} and  $|z|\ge n$ imply $|\Theta_t(z)|\ge |z|- |u(z)|\ge n-\vv,$  we derive 
\beq\label{XXF}  \P(\tau_n^x< t) \le  \ff{ V(x+u_0(x)) +c_3 t}{\inf_{|y|\ge n-\vv} H(V(y))}=:\vv_{t,n}(x),\ \ t\in [0,T].\end{equation} 
Since $\lim_{|x|\to\infty} H(V)(x)=\infty,$ we have $\lim_{n\to\infty} \vv_{t, n}(x)=0$.
 Therefore, $\tau_n^x\to T$ as $n\to\infty$ as desired.

(b) Let $X_t^x$ and $X_t^y$ solve \eqref{E1} with initial values $x,y$ respectively. Then 
$$X_t^{x,n}= X_t^x,\ \ X_t^{y,n}=X_t^y,\ \ t\in [0,T\land \tau_{n}^x \land\tau^y_{n}].$$
Combining  this with \eqref{C1} and \eqref{XXF} for some $c(n)>0$, we obtain
\beg{align*}
&\sup_{x,y\in B_k, |x-y|\le\vv} \E\Big[\sup_{t\in [0,T]} |X_t^x-X_t^y|\land 1\Big]\\
& \le \sup_{x,y\in B_k, |x-y|\le\vv} \Big\{\E\Big[\sup_{t\in [0,T]} |X_t^{x,n}-X_t^{y,n}|\land 1\Big]+ \P(\tau_{n}^x\land\tau_{n}^y<T)\Big\}\\
&\le c(n)\vv+  \vv_{T,n}(x) +\vv_{T,n}(y),\ \ n\ge 1.\end{align*}
By letting first $\vv\downarrow 0$ then $n\to\infty$, we prove assertion (1).

(c) Let $P_t^n$ be associated with $X_t^n$. 
 By  the Bismut formula in Theorem 1.1 (iii) of \cite{XXZZ},  
we find  some constant $c_n>0$ such that
$$\|\nn P_t^nf\|_\infty\le \ff{c_n}{\ss t} \|f\|_\infty,\ \  t\in (0,T], f\in \B_b(\R^d).$$ Equivalently, 
$$\|(P_t^n)^*\dd_x-(P_t^n)^*\dd_y\|_{var}\le \ff{c_n}{\ss t} |x-y|,\ \ x,y\in\R^d, t\in (0,T].$$ Next, 
by \eqref{XXF} and $X_t=X_t^n$ for $t\le \tau_n^x$, we obtain
$$|P_tf(x)-P_t^nf(x)|\le 2 \|f\|_\infty\P(\tau_n^x\le t)\le 2\|f\|_{\infty} \vv_{t,n}(x)\to 0\ \text{as}\ n\to\infty.$$
Then 
\beg{align*}&\limsup_{y\to x} \|P^*_t \dd_x  -P^*_t \dd_y\|_{var} \\
&\le\limsup_{n\to\infty} \limsup_{y\to x}\sup_{|f|\le 1} \Big\{ |P_t^n f(x)-P_t^n f(y)|+ 2 \|f\|_\infty \vv_{t,n}(x)+2 \|f\|_\infty \vv_{t,n}(y)\Big\}\\
&=0,\ \ t\in (0,T].\end{align*} 
So, assertion (2) is proved.

(d) When $\phi(r)=r$, by $(A_2)(3')$, \eqref{YNN0} and It\^o's formula,  for any $k\ge 1$ we find a constant $c_1(k)>0$ such that 
$$\d \{V(Y_t^n)^k\} \le c_1 (k)  V(Y_t^n)^k \d t + \d M_t^k$$
for some martingale $M_t^k$ with $\d\<M^k\>_t\le \{c_1(k) V(Y_t^n)^k\}^2\d t.$
Combining this with BDG's inequality, \eqref{*4} and $(A_2)(3)$, we find    constants $c_2(k), c_3(k) >0$ such that
\beg{align*} &\E\Big[\sup_{t\in [0,T]}V(X_t^n)^k\Big]\le (1+\vv K)  \E\Big[\sup_{t\in [0,T]}V(Y_t^n)^k\Big]\\
&\le c_2(k) V(x+ u_0^\ll(x))^k\le c_3(k)V(x)^k,\ \ n\ge 1\end{align*}
By Fatou's lemma with $n\to\infty$,  we prove  assertion (3) for some constant $c(k)>0$. 

\end{proof} 

\subsection{Wang's Harnack inequality} 

Under the monotone condition
 \beg{align*} &2\<b_t(x)-b_t(y), x-y\>+\|\si_t(x)-\si_t(y)\|_{HS}^2 \le K|x-y|^2,\\
 &|\{\si_t(x)-\si_t(y)\}^*(x-y)|\le K |x-y|,\ \ x,y\in \R^d, t\in [0,T],\end{align*}
 the following  Wang's Harnack   inequality was established in \cite{W11} for large $p>1$ and some constant $c>0$:
 $$ |P_t f(y)|^p\le \e^{\ff {c|x-y|^2}t} P_t|f|^p(x),\ \ x,y\in \R^d, t\in (0,T], f\in \B_b(\R^d).$$
Our next result   extends this inequality  to the singular setting, which generalizes   the main result in \cite{Shao} for Lipschitz continuous  $\si_t$ as well as  the corresponding result in
  \cite[Theorem 4.3(1)]{YZ0} for $\ff 1 2$-H\"older $\si_t$, since $(H2)(2)$ allows $\si_t$ to only have  a Dini type continuity.

\beg{thm}\label{T02} Assume $(A_1)$ and $(H2)(2).$  Then there exist constants $\hat p>1$ and $c>0$ such that 
   for all $p\ge \hat p,$
 \beq\label{HKK} |P_t f(y)|^p\le \e^{c+\ff {c|x-y|^2}t} P_t|f|^p(x),\ \ x,y\in \R^d, t\in (0,T], f\in \B_b(\R^d).\end{equation} 
\end{thm} 

\beg{proof} (a) We first observe that it suffices to prove for  $b^{(0)}=0$. Indeed, let $\hat P_t$ be the semigroup associated with the SDE
$$\d X_t^x = b_t^{(1)}(X_t^x)\d t+ \si_t(X_t^x)\d W_t, \ \ t\in [0,T].$$ Let
$$R^x:=\e^{\int_0^T \<\{\si_t^*(\si_t\si_t^*)^{-1}b_t^{(0)}\}(X_t), \d W_t\>-\ff 1 2 \int_0^T |\{\si_t^*(\si_t\si_t^*)^{-1}b_t^{(0)}\}(X_t)|^2\d t}.$$
By   $(A_1)$ and Khasminskii's estimate Lemma 4.1 in \cite{XXZZ}, we have 
$$\sup_{x\in \R^d} \E[ |R^x|^q]<\infty,\ \ q>1.$$
Then by Gisranov's theorem, for any $p>1$ there exists $c(p)>0$ such that 
$$|P_tf|^p(x)= \big|\E[R^x f(X_t^x)]\big|^p\le \big(\E \big[|R^x|^{\ff p{p-1}}\big]\big)^{p-1} \E [|f|^p(X_t^x)]\le c(p) \hat P_t |f|^p(x),\ \ p>1.$$
Similarly, the same inequality holds by exchanging positions of $P_t$ and $\hat P_t$. Thus, if the desired assertion holds for $\hat P_t$, it also holds for $P_t$.

(b) Now, we consider the regular case that   $b=b^{(1)}$. 
In this case,  there exists a constant $K>0$ such that for any $x,y\in\R^d,$
\beq\label{HA1} 2\<x-y,b_t(x)-b_t(y)\>+ \|\si_t(x)-\si_t(y)\|_{HS}^2 \le K (|x-y|^2\lor \varphi(|x-y|)^{2}).\end{equation}
For fixed $t\in (0,T]$, let 
\beq\label{GGS} \gg_s= \ff {1-\e^{K(s-t)}} K,\ \ s\in [0,t],\end{equation}  so that for some constant $K_1>1$ 
\beq\label{HA2} K\gg_s -2 -\gg_s'= -1,\ \ K_1(t-s) \ge \gg_s \ge K_1^{-1} (t-s),\ \ s\in [0,t].\end{equation} 
Since the coefficients of the following SDE are continuous and of linear growth in $x$ locally uniformly in $s\in [0,t)$, it has a weak solution (note that in general it is not well-posed)
\beq\label{CPL} \beg{cases} \d X_s = b_s(X_s)\d s +\si_s(X_s)\d W_s,\ \ &X_0=x,\\
\d Y_s = \big\{b_s(Y_s)+ \si_s(Y_s) \xi_s\big\}\d s +\si_s(Y_s)\d W_s,\ \ &Y_0=y, s\in [0,t),\end{cases}\end{equation}
where 
\beq\label{XI'} \xi_s:= \ff{\{\si_s^*(\si_s\si_s^*)^{-1} \}(X_s)(X_s-Y_s)}{\gg_s},\ \ s\in [0,t].\end{equation} 
The   coupling \eqref{CPL} is modified from \cite{W11}, we will show that it implies $X_t=Y_t$ which is crucial to establish Wang's Harnack inequality. 
Let \beq\label{RNN} \beg{split} & \tau_n= \ff{nt}{n+1} \land \inf\big\{s\ge 0:  |X_s|\lor |Y_s|\ge n\big\},\\
 &R_r:= \e^{-\int_0^{r} \<\xi_s,\d W_s\> -\ff 1 2\int_0^{\tau_n} |\xi_s|^2\d s},\ \ r\in [0,t].\end{split}\end{equation}
 By Girsanov's theorem,
$$\tt W_s:= W_s + \int_0^{s \land\tau_n} \xi_r\d r,\ \ s\in [0,t]$$ is an $m$-dimensional Brownian motion under the probability $\Q_n:= R_{\tau_n}\P$.
So,  before time $\tau_n$, \eqref{CPL} is reformulated as
$$\beg{cases} \d X_s =\big\{ b_s(X_s)-\ff{X_s-Y_s}{\gg_s}\big\}\d s +\si_s(X_s)\d \tt W_s,\ \ &X_0=x,\\
\d Y_s =  b_s(Y_s) \d s +\si_s(Y_s)\d \tt W_s,\ \ &Y_0=y, s\in [0,\tau_n].\end{cases}$$
By \eqref{HA1} and It\^o's formula, we obtain  
\beg{align*}   \d |X_s-Y_s|^2 \le  &\Big\{K (|X_s-Y_s|^2   +\varphi( |X_s-Y_s|)^{2})
 -\ff{2|X_s-Y_s|^2 }{\gg_s}\Big\}\d s + \d M_s\end{align*}
 for $s\in [0,\tau_n]$ and the $\Q_n$-martingale
$$\d M_s:= 2 \<X_s-Y_s, (\si_s(X_s)- \si_s(Y_s))\d \tt W_s\>$$ satisfying
\beq\label{MTT} \d\<M\>_s \le K^2 |X_s-Y_s|^2 \d s ,\ \ s\in [0,\tau_n].\end{equation}
Combining this with   It\^o's formula, we obtain 
\beq\label{**0} \beg{split} &\d\Big\{ \ff{|X_s-Y_s|^2}{\gg_s}\Big\}\\
&\le \Big\{\ff{(K\gg_s-2-\gg_s') (|X_s-Y_s|^2 }{\gg_s^2}+\ff{K\varphi( |X_s-Y_s|)^2}{\gg_s} \Big\}\d s +\frac{\d M_s }{\gg_s}. \end{split}\end{equation}

On the other hand, we observe that
\beq\label{**1}\beg{split} & \ff{K\varphi( |X_s-Y_s|)^2}{\gg_s} -\ff{|X_s-Y_s|^2}{2\gg_s^2} \\
&\le \sup_{r>0} \Big\{\ff{K\varphi(r)^2}{\gg_s}-\ff{r^2}{2\gg_s^2}\Big\}
\le \ff{K(\varphi\circ \psi^{-1})^2(2K\gg_s)}{\gg_s}:=g_t(s).\end{split}\end{equation} 
Indeed, since $\psi(r):=\ff {r^2}{\varphi(r)^2}$ is increasing in $r$, for $r\ge \psi^{-1} (2K\gg_s)$ we have 
\beg{align*} &\ff{K\varphi(r)^2}{\gg_s}-\ff{r^2}{2\gg_s} = \ff{K\varphi(r)^2}{\gg_s}\Big( 1-\ff{\psi(r)}{2K\gg_s}\Big) \\
&\le \ff{K\varphi(r)^2}{\gg_s}\Big( 1-\ff{\psi(\psi^{-1}(2K\gg_s))}{2K\gg_s}\Big) = 0,\end{align*}
while for $r<\psi^{-1} (2K\gg_s)$  
$$\ff{K\varphi(r)^2}{\gg_s}-\ff{r^2}{2\gg_s} \le  \ff{K\varphi(\psi^{-1}(2K\gg_s))^2}{\gg_s},$$
so that \eqref{**1} holds. Combining \eqref{**0} and \eqref{**1}, and noting that $\int_0^1\ff{\varphi(\psi^{-1}(s))^2}s\d s<\infty$ implies 
$$c_1:= \sup_{t\in [0,T]} \int_0^t g_t(s) \d s<\infty$$ 
by \eqref{HA2}, we have
$$\int_0^{\tau_n} \ff{|X_s-Y_s|^2}{2 \gg_s^2}\d s \le  \ff{|x-y|^2}{2 \gg_0^2}+c_1 +\int_0^{\tau_n} \ff{\d M_s}{\gg_s},\ \ s\in [0,\tau_n].$$
By this and \eqref{MTT},   for any $\ll>0$ we have 
\beg{align*}& \e^{-\Big(\ll c_1+\ff{\ll |x-y|^2}{\gg_0}\Big)} \E_{\Q_n} \Big[\e^{\ll \int_0^{\tau_n}  \ff{|X_s-Y_s|^2 }{\gg_s^2}\d s} \Big] \le \E_{\Q_n}\big[\e^{\ll \int_0^{\tau_n} \ff{\d M_s}{\gg_s}}\big]\\
&\le \big(\E_{\Q_n} \big[\e^{2\<M\>_{\tau_n}}\big]\big)^{\ff 1 2} \le \Big(\E_{\Q_n} \Big[\e^{2K^2\ll^2\int_0^{\tau_n}  \ff{|X_s-Y_s|^2}{\gg_s^2}\d s} \Big]\Big)^{\ff 1 2}.\end{align*}
Taking $\ll= (2K^2)^{-1}$ and noting that  \eqref{HA2} implies $\gg_0\ge K_1 t$, we find a constant $c_2>0$ such that  
\beq\label{WRE} \sup_{n\ge 1} \E_{\Q_n}  \Big[\e^{\ll \int_0^{\tau_n}  \ff{|X_s-Y_s|^2}{\gg_s^2}\d s} \Big]\le \e^{c_2+ \ff{c|x-y|^2}t}.\end{equation} 
Since \eqref{XI'} implies
$$|\xi_s|^2\le \ff{c_3|X_s-Y_s|^2}{\gg_s^2}$$
for some constant $c_3>0,$ this implies that for some constants $q,c_4>1,$
\beq\label{PAY} \sup_{n\ge 1} \E \big[|R_{\tau_n}|^q\big] \le \e^{c_4+ \ff{c_4|x-y|^2}t}.\end{equation} 
By the martingale convergence theorem, this implies that  $(R_s)_{s\in [0,t]}$ is a martingale with
\beq\label{RTT} \E[ R_t^q]\le  \e^{c_4+\ff{c_4|x-y|^2}t},\end{equation}
 such that Girsanov's theorem implies that $(\tt W_s)_{s\in [0,t]}$ is an $m$-dimensional Brownian motion under $\Q:=R_t\P$,   and
$$ \d Y_s =  b_s(Y_s) \d s +\si_s(Y_s)\d \tt W_s,\ \ Y_0=y, s\in [0,t] $$ holds 
so that $P_t f(y)= \E_\Q[f(Y_t)]$, and furthermore \eqref{PAY} ensures 
$$\E_\Q \Big[\e^{\ll \int_0^{t}  \ff{|X_s-Y_s|^2}{\gg_s^2}\d s} \Big]<\infty.$$
Since $\int_0^t \ff{\d s}{\gg_s^2}=\infty$ and $  |X_s-Y_s|^2 $ is continuous in $s$, this implies  $\Q(X_t=Y_t)=1$. 
Combining this with \eqref{RTT} and   $P_t f(y)= \E_\Q[f(Y_t)]$, we find a constant $c>0$ such that  for any $p\ge  \ff q {q-1}$, H\"older's inequality yields  
\beg{align*} &|P_t f(y)|^p= |\E[R_t f(Y_t)]|^p= |\E[R_tf(X_t)]|^p\\
&\le (\E R_t^q)^{\ff p q}  \E[ |f|^p(X_t)]\le (P_t|f|^p)(x) \e^{c+ \ff{c|x-y|^2}t}. \end{align*}
 \end{proof} 

\section{Proof  of Theorem \ref{TA1}  }

  \paragraph{ Proof of Theorem \ref{TA1}(1).}
  Let $X_0$ be $\F_0$-measurable with $\gg:=\L_{X_0}\in \scr P_V$. Let 
 $$\C^\gg:=\{\mu\in C([0,T];\scr P_V):\ \mu_0=\gg\}.$$
 For any $\mu\in C([0,T]; \scr P_V)$, by Theorem \ref{T01}, $(A_2)$   implies that the following SDE is well-posed
\beq\label{EM} \d X_t^\mu= b_t(X_t^\mu,\mu_t)\d t +\si_t(X_t^\mu)\d W_t,\ \ X_0^\mu=X_0.\end{equation} 
Denote $\Phi_t(\mu):= \L_{X_t^\mu}.$    By Theorem \ref{T01}, for the well-posedness of \eqref{E10} and estimate \eqref{KK2},   it suffices to prove that $\Phi$ has a unique fixed point in
 $\C^\gg$. To this end, following the line of \cite{HW20c} and \cite{W21b}, we approximate $\C_T^\gg$ by bounded subsets
 $$\C_N^\gg:=\Big \{\mu\in C^\gg: \sup_{t\in [0,T]} \mu_t(V) \e^{-Nt} \le N(1+\gg(V))\Big\},\ \ N\ge 1.$$

 (1a) We claimed that for some constant $N_0\ge 1$, $\Phi \C_N^T\subset \C_N^T$ for $N\ge N_0.$
 To this end,   let
 $$L_t^{\mu}:=\nn_{ b_t^{(0)}(\cdot,\mu_t)} +\ff 1 2 {\rm tr}\big\{ \si_t\si_t^* \nn^2\big\} $$
 and  consider the Zvonkin's transform of $X_t^{\mu}$ and the Kolmogorov backward equation as follows, 
\beq\label{ZOVKI}\beg{split} & Y_t^{\mu}=X_t^{\mu}+u_t^{\mu}(X_t^{\mu}),\ u_t^{\mu}\in \tilde{H}_{q_0}^{p_0,}\\
& (\pp_t+L_t^{\mu})u_t^{\mu}=\lambda u_t^{\mu}-b^{(0)}_t, \ \  t\in[0,T],\ u^{\mu}_T=0,
 \end{split}\end{equation}
for $\lambda>0$ such that $\|u_t^{\mu}\|_{\infty}+\|\nn u_t^{\mu}\|_{\infty}\leq \ff 1 2$.  
  By $(H_1)(3)$ and It\^o's formula, we find a constant $c_1>0$ such that
\beq\label{CD0} \d \{V(Y_t^\mu)\}^2 \le c_1 \{V(Y_t^\mu)^2+\mu_t(V)^2 \} \d t + \d M_t\end{equation} 
 for some martingale $M_t$. By the condition on $V$ and $|X_t^\mu-Y_t^\mu|\le \ff 1 2$, we find a constant $C>1$ such that
 \beq\label{NBB} C^{-1} V(X_t^\mu)\le V(Y_t^\mu)\le C V(X_t^\mu),\end{equation}
 so that \eqref{CD0}   implies that for some constant $c_2>0$ 
 \beq\label{CDS}\beg{split} & \E\big(V(X_t^\mu)^2\big|X_0^\mu\big)\le  C^2 \e^{c_1t} V(X_0^\mu)^2+ C^2c_1\int_0^t \e^{c_1(t-s)} \mu_s(V)^2 \d s \\
& \le c_2 V(X_0^{\mu})^2+ c_2\e^{c_2t} \{N(1+\gg(V)\}^2\int_0^t \e^{(2N-c_1)s}\d s\\
 &\le c_2 V(X_0^{\mu})^2+\frac {c_2}{(1-K\ee)^2}\e^{c_2t}  \{N(1+\gg(V)\}^2\ff{1}{2N-c_1} \e^{(2N-c_1)t} ,\ \ t\in [0,T], \ \ \mu\in \C_N^\gg.\end{split}\end{equation}
    Thus,  for  any $N\ge N_0:=c_2+2\ss{c_2},$  we have 
 $$\sup_{t\in [0,T]} \{\Phi_t(\mu)\}(V)\e^{- Nt} \le   (1+\gg(V))\big\{\ss{c_2}  + \ss {c_2 N}\big\} \le N(1+\gg(V)).$$ 
 Thus, $\Phi \C_N^\gg \subset \C_N^\gg$ for $N\ge N_0.$ 

 (b) Let $N\ge N_0$. We prove that $\Phi$ has a unique fixed point in $\C_N^\gg$, and hence it has a unique fixed point in $\C^\gg$ as desired. 
 Consider the following  complete metric on $\C_N^\gg$:
 $$\rr_{\ll}(\mu,\nu):= \sup_{t\in [0,T]} \e^{-\ll t} \|\mu_t-\nu_t\|_V.$$
 Let  $$ \xi_s:= \{\si_s^*(\si_s\si_s^*)[b_s(X_s^\mu, \nu_s)- b_s(X_s^\mu, \mu_s)]\}(X_s^\mu),\ \ s\in [0,T].$$ 
 By $\eqref{KK1}$, 
\beq\label{NRT} R_t:=\e^{\int_0^T \<\xi_s, \d W_s\>-\ff 1 2 \int_0^T |\xi_s|^2\d s}\end{equation} 
 is a martingale,  such that 
 $$\tt W_r:= W_r-\int_0^r \xi_s \d s,\ \ r\in [0,t]$$
 is a Brownian motion under the probability $\Q_t:= R_t\P.$ 
 Reformulate \eqref{EM} as 
 $$\d X_r^\mu= b_r(X_r^\mu,\nu_r)\d r +\si_r(X_r^\mu)\d \tt W_r,\ \ X_0^\mu=X_0,\ \ r\in [0,t].$$
   By the uniqueness we obtain
   $$\Phi_t(\nu)= \L_{X_t^\nu}= \L_{X_t^\mu|\Q_t},$$ where $\L_{X_t^\mu|\Q_t}$ stands for the distribution of $X_t^\mu$ under $\Q_t$. 
   Then by \eqref{CDS}, we find a constant $c_1(N)>0$ 
\beq\label{RPP2} \beg{split} &\|\Phi_t(\mu)-\Phi_t(\nu)\|_V = \sup_{|f|\le V} \big|\E\big[f(X_t^\mu)(1-R_t)\big]\big| \\
&\le \E\Big[\big\{\E(V(X_t^\mu)^2|X_0^\mu)\big\}^{\ff 1 2} \big\{\E [|R_t-1|^2|X_0^\mu]\big\}^{\ff 1 2}\Big]\\
   &\le c_1(N) \E \Big[V(X_0) \big\{\E [R_t^2-1|X_0]\big\}^{\ff 1 2}\Big].\end{split} \end{equation} 
Since $\mu\in \C_N^\gg$, by $\eqref{KK1}$ we 
find a constant $c_2(N)>0$ such that 
$$|\xi_s|^2\le c_2(N) (1\land \|\mu_s-\nu_s\|_V^2),\ \ s\in [0,T],$$
so that  for some constant $c_3(N)>0$ 
\beg{align*} &\E [R^2_t-1|X_0]\le \E \bigg[\e^{2 \int_0^t \<\xi_s,\d W_s\>- \int_0^t |\xi_s|^2\d s} -1\bigg|X_0\bigg]\\
&\le \E \bigg[\e^{2 \int_0^t \<\xi_s,\d W_s\>-2 \int_0^t |\xi_s|^2\d s+\int_0^t |\xi_s|^2\d s } \bigg]\le \E \big[\e^{2 \int_0^t |\xi_s|^2\d s}\big|X_0\big] -1\\
&\le c_3(N) \int_0^t \|\mu_s-\nu_s\|_V^2\d s .\end{align*}   
Combining this with \eqref{RPP2}, we find a constant $c_4(N)>0$ such that 
\beg{align*} &\rr_\ll(\Phi(\mu), \Phi(\nu))=\sup_{t\in [0,T]} \e^{-\ll t} \|\Phi_t(\mu)-\Phi_t(\nu)\|_V\\
&\le c_4(N) \E[V(X_0)] \rr_\ll(\mu,\nu)\sup_{t\in [0,T]}\bigg(\int_0^t \e^{-2\ll(t-s)}\d s\bigg)^{\ff 1 2} \\
&\le \ff{c_4(N) \E[V(X_0)] }{\ss{2\ll}} \rr_\ll(\mu,\nu),\ \ \mu,\nu\in \C_N^\gg.\end{align*} 
 Therefore, when $\ll>0$ is large enough, $\Phi$ is contractive under $\rr_\ll$ so that it has a unique fixed point in $\C_N^\gg$ as desired. 
 
 (c) Proof of \eqref{KK2}. Let $X_t$ solve \eqref{E10} with $\L_{X_0}\in \scr P_V$, and denote $\mu_t=\L_{X_t}.$ We have $\sup_{t\in [0,T]}\mu_t(V)<\infty.$ By $(H_1)(3)$ and It\^o's formula, we find a constant $c_1>0$ such that 
 \beq\label{ER1}  \d V(Y_t^{\mu})\le c_1 \big\{V(Y_t^{\mu})+ \mu_t(V)\big\}\d t + \d M_t\end{equation} for some martingale $M_t$ with
 \beq\label{ER*} \d\<M\>_t\le c_1^2 V(X_t)^2 \d t.\end{equation}
 By this and \eqref{NBB}, we find a constant $c_2>0$ such that  
 \beq\label{VES} \mu_t(V)= \E[V(X_t)]\le  c_2 \int_0^t \mu_s(V)\d s,\ \ t\in [0,T],\end{equation}
  so that by Gronwall's inequality, 
 \beq\label{ERR} \E[V(X_t)]\le \e^{c_2 t} \E[V(X_0)],\ \ t\in [0,T].\end{equation}
Combining this with   \eqref{ER1} and applying It\^o's formula, for any $p\ge 1$ we find a constant $c_1(p)>0$ such that 
 $$\d V(Y_t^{\mu})^p \le c_1(p)  \big\{V(Y_t^{\mu})^p+ \mu_t(V)^p\big\}\d t +pV(Y_t^{\mu})^{p-1} \d M_t.$$
By \eqref{ER*}, \eqref{NBB} -and   BDG's inequality, we find a constant $c_2(p)>0$ such that 
 $$\xi_t:= E\bigg[\sup_{s\in [0,t]} V(X_s)^p\Big|X_0\bigg],\ \ t\in [0,T]$$ satisfies 
\beg{align*}  \xi_{t\land\tau_n} & \le V(X_0)^p+ c_2(p)(\E [V(X_0)] )^p + c_2(p)  \int_0^t \xi_{s\land\tau_n}\d s + c_2(p) \E\bigg[\bigg(\int_0^{t \land \tau_n} V(X_s)^{2p} \d s \bigg)^{\ff 1 2}\bigg|X_0\bigg]\\
 &\le  V(X_0)^p+c_2(p)\E [V(X_0)] + c_2(p)  \int_0^t \xi_{s\land\tau_n}\d s + \ff 1 2 \xi_{t\land\tau_n} +\ff { c_2(p)^2} 2  \int_0^t \xi_{s\land \tau_n}\d s,\ \ t\in [0,T].\end{align*}
 So that for $c_3(p):= 2 c_2(p)+c_2^2(p)$ we obtain
 $$\xi_{t\land\tau_n}\le 2\big\{V(X_0)^p+  c_2(p)(\E[V(X_0)])^p\big\}   \e^{c_3(p) t},\ \ t\in [0,T],n\ge 1.$$
 Letting $n\to\infty$ we derive \eqref{KK2} for some constant $c(p)>0$. 
 
 \paragraph{ Proof of Theorem \ref{TA1}(2).}  Let $\hat P_t$ be the Markov semigroup of $X_t^\mu$ solving \eqref{EM} for $\mu_t:= P_t^*\mu$, so that 
 \beq\label{PPY}  P_t^*\mu= \hat P_t^*\mu,\ \ t\in [0,T].\end{equation} 
 By Theorem \ref{T01}, we have 
 $$\lim_{y\to x} \|\hat P_t^*\dd_x-\hat P_t^*\dd_y\|_{var}=0,\ \ x\in \R^d.$$
Since $\mu_n\to \mu$ weakly, we may construct random variables $\{\xi_n\}$ and $\xi$ such that $\L_{\xi_n}=\mu_n, \L_{\xi}=\mu$ and 
$\xi_n\to\xi$ a.s. Thus, by the dominated convergence theorem we obtain
\beq\label{RPP0}  \beg{split}&   \lim_{n\to\infty} \|\hat P_t^*\mu_n- \hat P_t^*\mu\|_{var} = \lim_{n\to\infty} \|\E[\hat P_t^*\dd_{\xi_n} - \hat P_t^*\dd_{\xi}]\|_{var} \\
&\le  \lim_{n\to\infty} \E[\|\hat P_t^*\dd_{\xi_n} - \hat P_t^*\dd_{\xi}\|_{var}]=0.\end{split} \end{equation} 
Hence,
\beq\label{RPP3}\beg{split}& \limsup_{n\to\infty}  \|\hat P_t^*\mu_n  -\hat P_t^*\mu \|_{V}\\
&\le \limsup_{n\to\infty}\Big\{ \sup_{|f|\le N}   |(\hat P_t^*\mu_n)(f)  -(\hat P_t^*\mu)(f) |+ 
   \|\hat P_t^*\mu_n  -\hat P_t^*\mu \|_{V}\Big\}\\
&\le N  \limsup_{n\to\infty} \|\hat P_t^*\mu_n  -\hat P_t^*\mu\|_{var} + \sup_{n\ge 1} \int_{\R^d} \hat P_t (V-N)^+ \d (\mu_n+\mu) \\
&=  \sup_{n\ge 1}    \big\{\hat P_t^*(\mu_n +\mu)\big\}\big((V-N)^+ \big),\ \ N\ge 1.\end{split}\end{equation} 
Since $\mu_n(V^p)$ is bounded for some $p\in (1,2]$,  \eqref{CDS} implies that
\beq\label{RPPM} \sup_{n\ge 1, t\in [0,T]} (\hat P_t^* \mu_n)(V^p)<\infty,\end{equation}  so that 
letting $m\to\infty$ in \eqref{RPP3} we prove  
\beq\label{RPPN} \limsup_{n\to\infty}  \|\hat P_t^*\mu_n  -\hat P_t^*\mu \|_{V}=0.\end{equation}

On the other hand, by the Girsanov transform in step (b) above for $\mu_n$ replacing $\nu$, we find a constant $c>0$ such that 
$$\|P_t^*\mu_n-\hat P_t^*\mu_n\|_{V}^2\le c \int_0^t \|\mu_s-\nu_s\|_V^2\d s,\ \ t\in [0,T].$$
Combining this with \eqref{RPPN} and Fatou's lemma due to \eqref{RPPM}, we derive 
\beg{align*} &\limsup_{n\to\infty} \|P_t^*\mu_n-\hat P_t^*\mu\|_{V}^2\le 2  \limsup_{n\to\infty}  \big\{\|\hat P_t^*\mu_n  -\hat P_t^*\mu \|_{V}^2+  \|P_t^*\mu_n-\hat P_t^*\mu_n\|_{V}^2\big\}\\
&\le 2\int_0^t \limsup_{n\to\infty} \|P_s^*\mu_n-\hat P_s^*\mu_n\|_{V}^2\d s<\infty,\ \ t\in [0,T],\end{align*} 
By Gronwall's inequality and \eqref{PPY}, we obtain 
$$\limsup_{n\to\infty} \|P_t^*\mu_n- P_t^*\mu\|_{V}=\limsup_{n\to\infty} \|P_t^*\mu_n-\hat P_t^*\mu\|_{V}=0.$$
This implies  \eqref{KK3}.

\paragraph{ Proof of Theorem \ref{TA1}(3).}  By \eqref{KK1'},     $R_t$ in \eqref{NRT}  is a martingale with
 $$|\xi_s|^2\le c \|\mu_s-\nu_s\|_{var}^2$$
 for some constant $c>0$. Then 
 \beq\label{RPP4} \|\Phi_t(\mu)-\Phi_t(\nu)\|_{var}=\sup_{|f|\le 1} |\E[f(X_t^\mu)(R_t-1)]|\le \E[|R_t-1|].\end{equation} 
 By Pinsker's inequality, we obtain 
 $$(\E[|R_t-1|])^2\le 2 \E[R_t\log R_t] = 2\E_{\Q_t } \log R_t]= \E_{\Q_t} \int_0^t |\xi_s|^2 \d s\le c^2 \int_0^t \|\mu_s-\nu_s\|_{var}^2\d s.$$
 Combining this with \eqref{RPP4}, as shown in the proof of (1) we see that when $\ll>0$ is large enough, $\Phi$ is contractive in $\C^\gg$ under the metric
 $$\tt\rr_{\ll}(\mu,\nu):=\sup_{t\in [0,T]} \e^{-\ll t} \|\mu_t-\nu_t\|_{var}.$$ 
 Hence, by Theorem \ref{T01}, \eqref{E10} is well-posed  and \eqref{KK2} holds. 
 
Let $\mu_t=P_t^*\mu$ as in the proof of  Theorem \ref{TA1}(3).  As shown above  that \eqref{KK1'}, Girsanov's theorem and Pinsker's  inequality 
 imply
 $$\|P_t^*\mu_n-\hat P_t^*\mu_n\|_{var}^2\le c \int_0^t \|P_s^*\mu_n- P_s^*\mu\|_{var}^2\d s$$ for some constant $c>0$. 
Thus, by the same reason leading to \eqref{KK3},  \eqref{KK3'} follows from \eqref{RPP0}.

\section{Proof of Theorem \ref{TA2}}

Noting that conditions in Theorem \ref{TA2} imply those in Theorem \ref{TA1} and when $V$ is bounded we have 
$$\|\cdot\|_{var}\le \|\cdot\|_V\le \|V\|_\infty \|\cdot\|_{var},$$
so the first assertion follows. It remains to verify  \eqref{HHH} and \eqref{GRR}.

(1) Let $\hat P_t$ be associated to solutions of \eqref{E1} for $b(\cdot,\dd_0)$ replacing $b$. By Theorem  \ref{T02},
there exist constants $c',p'>1$ such that 
$$ |\hat P_t f(y)|^{p'}\le \e^{c'+\ff {c'|x-y|^2}t} \hat P_t|f|^{p'}(x),\ \ x,y\in \R^d, t\in (0,T], f\in \B_b(\R^d).$$
Consequently, 
\beq\label{LNM} |\hat P_t f(\mu)|^{2p'}\le C(t,\mu,\nu)  \hat P_t|f|^{2p'}(\nu),\ \ x,y\in \R^d, t\in (0,T], f\in \B_b(\R^d)\end{equation} holds for
$$C(t,\mu,\nu):= \inf_{\pi\in \C(\mu,\nu)}\int_{\R^d\times \R^d} \e^{2c'+\ff {c'|x-y|^2}t}\pi(\d x,\d y).$$

Next, let $\hat X_t$ solve \eqref{E1} for $b(\cdot,\dd_0)$ replacing $b$ with initial distribution $\gg\in \scr P_V$, and denote
$$\xi_t:=\big\{\si_t^*(\si_t\si_t^*)^{-1}[b_t(\hat X_t, P_t^*\gg)-b_t(\hat X_t,\dd_0)]\big\},\ \ t\in [0,T].$$
By \eqref{MOO},
$$R_t:=\e^{\int_0^t\<\xi_s,d W_s\>-\ff 1 2 \int_0^t |\xi_s|^2\d s},\ \ t\in [0,T]$$ is a martingale, and by Girsanov's theorem,  
  for any $q>1$ we find a constant $c(q)>0$ such that 
$$|P_t f(\gg)|^q= |\E [R_t f(\hat X_t)]|^q\le \hat  P_t|f|^q(\gg) \e^{c(q) (\gg^2(V)+V(0))t},\ \ t\in [0,T].$$
Similarly, 
$$|\hat P_t f(\gg)|^q \le  P_t|f|^q(\gg) \e^{c(q) (\gg^2(V)+V(0))t},\ \ t\in [0,T].$$
Combining these with \eqref{LNM}, we derive \eqref{HHH} for any $p>p'$ and some constant $c>0$. 

(2) When $\Phi$ is bounded,  \eqref{MOO} implies 
$$\sup_{\mu\in \scr P} \|b(\cdot,\mu)-b(\cdot,\dd_0)\|_\infty<\infty.$$ Let $P_t^\mu$ be the Markov semigroup for solutions to \eqref{E1} for $b_t(\cdot,P_t^*\mu)$ replacing $b_t$,
 by \cite[Theorem 4.1]{YZ0}, there exists a constant $c'>0$ such that $P_t^\mu$ satisfies the log-Harnack inequality 
$$ P_t^\mu \log f(x)\le \log  P_t^\mu f(y)+\ff {c'|x-y|^2}t,\ \ x,y\in \R^d, t\in (0,T], 0<f\in \B_b(\R^d).$$
Consequently, 
$$  P_t^\mu \log f(\mu)\le \log  P_t^\mu f(\nu)+\ff {c'\W_2(\mu,\nu)^2}t,\ \ \mu,\nu\in \scr P, t\in (0,T], 0<f\in \B_b(\R^d).$$
Since $(P_t^\mu)^*\mu= P_t^*\mu$, this and Pinsker's inequality imply
$$ \|P_t^*\mu-(P_t^\mu)^*\nu\|_{var}^2 \le 2 \Ent((P_t^\mu)^*\mu|(P_t^\mu)^*\nu)=2\sup_{ P_t^\mu f(\nu)\le 1}  P_t^\mu \log f(\mu) \le \ff {2c'\W_2(\mu,\nu)^2}t.$$
Since $\|\cdot \|_{var}\le 2$, this is equivalent to
\beq\label{POO} \|P_t^*\mu-(P_t^\mu)^*\nu\|_{var}^2\le \aa_t:= \min\Big\{4, \ff {2c'\W_2(\mu,\nu)^2}t\Big\},\ \ t\in (0,T].\end{equation} 

On the other hand, let $X_t^{\mu,\nu}$ solve \eqref{E1} with $b_t(\cdot\mu_t)$ replacing $b_t$ and $\L_{X_0^{\mu,\nu}}=\nu$. Let
$$\xi_t:=\big\{\si_t^*(\si_t\si_t^*)^{-1}[b_t(X_t^{\mu,\nu}, P_t^*\nu)-b_t(X_t^{\mu,\nu},P_t^*\mu)]\big\},\ \ t\in [0,T].$$
By \eqref{MOO} for bounded $\Phi$,  we find a constant $K>0$ such that 
$$|\xi_t|^2\le K \|P_t^*\mu- P_t^*\nu\|_{var}^2.$$
So, $$R_t:=\e^{\int_0^t\<\xi_s,d W_s\>-\ff 1 2 \int_0^t |\xi_s|^2\d s},\ \ t\in [0,T]$$ is a martingale, and by Girsanov's theorem and Pinsker's inequality, we obtain 
$$\|(P_t^\mu)^* \nu- P_t^*\nu\|_{var}^2 \le 2 \E_{R_t\P}[\log R_t]\le  K \int_0^t \|P_s^*\mu-P_s^*\nu\|_{var}^2\d s.$$
Combining this with \eqref{POO}, and  we derive
\beg{align*} &\|P_t^* \mu- P_t^*\nu\|_{var}^2\le 2\|P_t^*\mu-(P_t^\mu)^*\nu\|_{var}^2+ 2\|(P_t^\mu)^* \nu- P_t^*\nu\|_{var}^2\\
&\le 2 \aa_t  +  2 K \int_0^t \|P_s^*\mu-P_s^*\nu\|_{var}^2\d s.\end{align*}
By Gronwall's inequality, we find a constant $c>0$ such that  
$$\|P_t^* \mu- P_t^*\nu\|_{var}^2\le \aa_t +2 K \int_0^t\aa_s \e^{2K(t-s)}\d s \le c\big(t^{-1}-\log [1\land \W_2(\mu,\nu)]\big) \W_2(\mu,\nu)^2.$$


\end{document}